\documentclass{elsart}

\usepackage{amsmath,amssymb,amsfonts,amsopn,amscd}
\usepackage{bbm,wasysym,stmaryrd}
\usepackage[english,francais]{babel}

\newtheorem{theorem}{Theorem}[section]
\newtheorem{lemma}[theorem]{Lemma}
\newtheorem{e-proposition}[theorem]{Proposition}
\newtheorem{corollary}[theorem]{Corollary}
\newtheorem{e-definition}[theorem]{Definition\rm}
\newtheorem{remark}{\it Remark\/}

\newtheorem{theoreme}{Th\'eor\`eme}[section]

\newtheorem{proposition}[theoreme]{Proposition}

\newtheorem{definition}[theoreme]{D\'efinition\rm}

\newcommand{\olivier}[1]{#1}
\newcommand{\of}[1]{#1}

\newcommand{\R}{\mathbbm{R}}
\newcommand{\mN}{\mathcal{N}}
\newcommand{\J}{J}
\newcommand{\Z}{\mathbbm{Z}}
\newcommand{\T}{\mathcal{T}}  
\newcommand{\mS}{\mathcal{S}}  
\newcommand{\M}{\mathcal{M}}

\newcommand{\mmeas}{\mathcal{M}_{\mS}}
\newcommand{\B}{\mathcal{B}}
\newcommand{\mM}{\mathcal{M}}

\newcommand{\undt}[1]{\underline{#1}}
\newcommand{\mone}{1}
\newcommand{\lsup}[1]{\underset{#1\to\infty}{\overline{\lim}}}

\def\og{\leavevmode\raise.3ex\hbox{$\scriptscriptstyle\langle\!\langle$~}}
\def\fg{\leavevmode\raise.3ex\hbox{~$\!\scriptscriptstyle\,\rangle\!\rangle$}}

\journal{the Acad\'emie des sciences}
\begin{document}
\centerline{}
\begin{frontmatter}


\selectlanguage{english}
\title{Asymptotic description of stochastic neural networks. II - Characterization of the limit law}


\selectlanguage{english}
\author[authorlabel1]{Olivier Faugeras},
\ead{firstname.name@inria.fr}
\author[authorlabel1]{James Maclaurin}
\address[authorlabel1]{Inria Sophia-Antipolis M\'editerran\'ee\\
NeuroMathComp Group}



\begin{abstract}
\selectlanguage{english}
We continue the development, started in \cite{faugeras-maclaurin:14}, of the asymptotic description of certain stochastic neural networks. We use the Large Deviation Principle (LDP) and the good rate function $H$ announced there to prove that $H$ has a unique minimum $\mu_e$, a stationary measure on the set of trajectories $\T^\Z$. We characterize this measure by its two marginals, at time 0, and from time 1 to $T$. The second marginal is a stationary Gaussian measure. With an eye on applications, we show that its mean and covariance operator can be inductively computed. Finally we use the LDP to establish various convergence results, averaged and quenched.
\vskip 0.5\baselineskip

\selectlanguage{francais}
\noindent{\bf R\'esum\'e} \vskip 0.5\baselineskip \noindent
\begin{center}{\bf  Description asymptotique de r\'eseaux de neurones stochastiques. II - caract\'erisation de la loi limite }\end{center}
Nous prolongeons le d\'eveloppement, commenc\'e en \cite{faugeras-maclaurin:14}, de la description asymptotique de certains r\'eseaux de neurones stochastiques. Nous utilisons le Principe de Grandes D\'eviations (PGD) et la bonne fonction de taux $H$ que nous y annoncions pour d\'emontrer l'existence d'un unique minimimum, $\mu_e$ de $H$, une mesure stationnaire sur l'ensemble $\T^\Z$ des trajectoires. Nous caract\'erisons cette mesure par ses deux maginales, \`a l'instant 0, et du temps 1 au temps $T$. La seconde marginale est une mesure gaussienne stationnaire. Avec un oeil sur les applications, nous montrons comment calculer de mani\`ere inductive sa moyenne et son op\'erateur de covariance. Nous montrons aussi comment utiliser le PGD pour \'etablir des r\'esultats de convergence en moyenne et presque s\^urement.
\end{abstract}
\end{frontmatter}

\selectlanguage{francais}
\section*{Version fran\c{c}aise abr\'eg\'ee}
Apr\`es avoir rappel\'e dans la section \ref{sec:maths} les notations et le mod\`ele de r\'eseaux de neurones utilis\'es dans \cite{faugeras-maclaurin:14}, nous montrons dans la proposition \ref{prop:unique} et le th\'eor\`eme \ref{theo:effective} que la bonne fonction de taux $H$ du PGD annonc\'e dans cette publication admet un minimum unique. Le th\'eor\`eme \ref{theo:effective} fournit une m\'ethode constructive de calcul effectif de la loi de ce minimum. Nous montrons enfin, dans la section \ref{sec:convergence}, l'int\'er\^et de ce minimum qui appara\^it comme la limite faible quand $n \to \infty$ (le nombre de neurones tend vers l'infini) de la loi $Q^{V_n}$ du r\'eseau moyenn\'ee par rapport aux poids synaptiques, c'est un r\'esultat en moyenne. Nous montrons aussi dans le corollaire \ref{cor:weakconv} un r\'esultat de convergence faible presque s\^urement par rapport aux poids synaptiques, r\'esultat int\'eressant d'un point de vue pratique puisqu'il \'evite de prendre la moyenne par rapport \`a tous les r\'eseaux. Le th\'eor\`eme \ref{theo:ergodic} donne un r\'esultat de convergence presque s\^ure de la mesure empirique vers le minimum de $H$.
\selectlanguage{english}
\section{Introduction}
In \cite{faugeras-maclaurin:14} we started our asymptotic analysis of very large networks of neurons with correlated synaptic weights. We showed that the image $\Pi^n$ of the averaged law $Q^{V_n}$ through the empirical measure satisfied a large deviation principle with good rate function $H$. In the same article we provided an analytical expression of this rate function in terms of the spectral representation of certain Gaussian processes. In the next section we recall some definitions given in  \cite{faugeras-maclaurin:14}.\vspace{-0.5cm}
\section{Mathematical framework}\label{sec:maths}
For some topological space $\Omega$ equipped with its Borelian $\sigma$-algebra $\B(\Omega)$, we denote the set of all probability measures by $\M(\Omega)$. We equip $\M(\Omega)$ with the topology of weak convergence. For some positive integer $n> 0$, we let $V_n = \lbrace j\in\Z: |j| \leq n \rbrace$, and $|V_n|=2n+1$. Let $\T = \R^{T+1}$, for some positive integer $T$. We equip $\T$ with the Euclidean topology, $\T^{\Z}$ with the cylindrical topology, and denote the Borelian $\sigma$-algebra generated by this topology by $\B(\T^{\Z})$. For some $\mu\in\M(\T^{\Z})$ governing a process $(X^j)_{j\in\Z}$,  we let $\mu^{V_n}\in \M(\T^{V_n})$ denote the marginal governing $(X^j)_{j\in V_n}$. For some $X \in \T$ and $0 \leq a \leq b \leq T$, $X_{a,b}$ denotes the $b-a+1$-dimensional subvector of $X$. We let $\mu_{a,b} \in \M(\T^\Z_{a,b})$ denote the marginal governing $(X^j_{a,b})_{j \in \Z}$. For some $j\in \Z$, let the shift operator $\mS^j:\T^{\Z}\to \T^{\Z}$ be $S(\omega)^k = \omega^{j+k}$. We let $\M_\mS(\T^{\Z})$ be the set of all stationary probability measures $\mu$ on $(\T^{\Z},\B(\T^{\Z}))$ such that for all $j\in\Z$, $\mu\circ (\mS^j)^{-1} = \mu$.  
Let $p_n:\T^{V_n} \to \T^{\Z}$ be such that $p_n(\omega)^k = \omega^{k\mod V_n}$. Here, and throughout the paper, we take $k\mod V_n$ to be the element $l\in V_n$ such that  $l = k \mod |V_n|$.  Define the process-level empirical measure $\hat{\mu}_n: \T^{V_n} \to \M_\mS\left(\T^{\Z}\right)$ as
\begin{equation}
\hat{\mu}_n(\omega) = \frac{1}{|V_n|}\sum_{k\in V_n} \delta_{S^k p_n(\omega)}.\label{defn:hatmun}
\end{equation}

The equation describing the time variation of the membrane potential
$U^j$ of the $j$th neuron writes
\begin{equation}\label{eq:U}
U^j_{t}=\gamma U^j_{t-1}+\sum_{i\in V_n} J_{ji}^{n}
f(U^i_{t-1})+\theta^j+B^j_{t-1}, , \quad U^j_0=u^j_0, \quad j \in V_n \quad t=1,\ldots,T
\end{equation}
\vspace{-0.5cm}

$f: \R \to ]0,\,1[$ is a monotonically increasing Lipschitz continuous bijection.
$\gamma$ is in $[0,1)$ and determines the time
scale of the intrinsic dynamics of the
neurons.
The $B^j_t$s are i.i.d. Gaussian random variables distributed as
$\mathcal{N}_1(0,\sigma^2)$\footnote{We note $\mN_p(m,\Sigma)$ the \of{density} of the $p$-dimensional
Gaussian variable with mean $m$ and covariance matrix $\Sigma$.}. They represent the fluctuations of the neurons' membrane potentials. The $\theta^j$s are i.i.d. as $\mathcal{N}_1(\bar{\theta},\theta^2)$. The are independent of the $B^i_t$s and represent the current injected in the neurons. The $u^j_0$s are i.i.d. random variables each governed by law $\mu_I$.

The $J_{ij}^{n}$s are the synaptic weights. $J_{ij}^{n}$ represents the
strength with which the `presynaptic' neuron $j$ influences the
`postsynaptic' neuron $i$. They arise from a stationary Gaussian random field specified by its mean \of{$\bar{J}$} and covariance function \of{$\Lambda : \Z\times \Z \to \R$}, see \cite{faugeras-maclaurin:14d,faugeras-maclaurin:14}. 
%

We note $\J^n$ the $|V_n| \times |V_n|$ matrix of the synaptic weights, $\J^n=(J_{ij}^n)_{i,j \in V_n}.$

where $\Psi:\T \rightarrow \T$ is the following affine bijection. Writing $v = \Psi(u)$, we define
\begin{equation}\label{eq:Psi}
\left\{
\begin{array}{lcl}
v_0 &=& \Psi_0(u)=u_0\\
v_s  &=&  \Psi_s(u)=u_s - \gamma u_{s-1}-\bar{\theta}\quad s=1,\cdots, T.
\end{array}
\right.
\end{equation}
For $v\in\T$, we write $\Psi^{-1}(v) = (\Psi^{-1}(v)_0,\ldots,\Psi^{-1}(v)_T)$. The coordinate $\Psi^{-1}(v)_{t}$ is the affine function of $v_s$, $s=0\cdots t$ obtained from equations 
\eqref{eq:Psi}
\begin{equation*}
\Psi^{-1}(v)_{t}=\sum_{i=0}^{t} \gamma^i v_{t-i}+\bar{\theta} \frac{\gamma^{t}-1}{\gamma-1}.
,\ t=0,\cdots,T.
\end{equation*}
We extend $\Psi$ to a mapping $\T^\Z \to \T^\Z$ componentwise and introduce the following notation.
\begin{definition}\label{def:mubarbar}
For each measure $\mu \in \M(\T^{V_n})$ or $\mmeas(\T^\Z)$ we define $\undt{\mu}$ to be $\mu \circ \Psi^{-1}$. 
\end{definition}
\of{Note that the correspondence $\mu \to \undt{\mu}$ is an isomorphism.}


We note $Q^{V_n}(\J^n)$ the element of $\M(\T^{V_n})$ which is the law of the solution to \eqref{eq:U} conditioned on $\J^n$. We let \of{$Q^{V_n}$} be the law averaged with respect to the weights. 

Finally we introduce the image law in terms of which the principal results of this paper are formulated.
\begin{definition}\label{def:PiNQN}
Let $\Pi^n \in \mathcal{M}(\mathcal{M}_\mS(\T^\Z))$ be the image law of $Q^{V_n}$ through the function
$\hat{\mu}_n: \T^{V_n} \to \mM_\mS(\T^\Z)$ defined by \eqref{defn:hatmun}:
\[
\Pi^n=Q^{V_n} \circ \hat{\mu}_n^{-1} 
\]
\end{definition}\vspace{-1.5cm}

\section{Characterization of the unique minimum of the rate function} 
In \cite{faugeras-maclaurin:14d}, to each measure $\nu \in \mM(\T^\Z)$ we associate the measure, noted $Q^\nu$ of  $\mM(\T^\Z)$ such that $\undt{Q}^\nu=\mu_I^\Z \otimes \undt{Q}^\nu_{1,T}$ where $\undt{Q}^\nu_{1,T}$ , a Gaussian measure on $\T_{1,T}^\Z$ with spectral density $\sigma^2 \delta(\theta)+\tilde{K}^\nu(\theta)$. The spectral density $\tilde{K}^\nu$ is defined in \cite{faugeras-maclaurin:14}. We also define the rate function $H^\nu$, which is a linear approximation of the functional $\Gamma$ defined in \cite{faugeras-maclaurin:14d} \of{and satisfies the relation $H^\mu(\mu)=H(\mu)$}. 
We prove the following lemma in \cite{faugeras-maclaurin:14d}.
\begin{lemma}
For $\mu,\nu\in\mM_\mS(\T^\Z)$, $H^\nu(\mu) = 0$ if and only if $\mu = Q^{\nu}$.
\end{lemma}
As stated in the following proposition, there exists a unique minimum $\mu_e$ of the rate function. We provide explicit equations for $\mu_e$ which would facilitate its numerical simulation.
\begin{proposition}\label{prop:unique}
There is a unique distribution $\mu_e \in \mathcal{M}_\mS(\T^\Z)$ which minimises $H$. This distribution satisfies $H(\mu_e)=0$ which is equivalent to
$\mu_e=Q^{\mu_e}$.
\end{proposition}
\begin{pf}
\of{The proof, which is found in \cite{faugeras-maclaurin:14d}, is an easy consequence of the explicit method we outline to actually calculate $\mu_e$ below in theorem \ref{theo:effective}.}
\qed
\end{pf}
We characterize the unique measure ${\mu_e}$ such that ${\mu_e}=Q^{\mu_e}$ in terms of its image $\undt{\mu_e}$. This characterization allows one to directly numerically calculate $\mu_e$. \of{Since $\undt{\mu_e}$ is Gaussian, the problem becomes that of defining the latest entries of $K^{\mu_e}$ and $c^{\mu_e}$ in terms of previous ones.} Hence we characterize $\undt{\mu_e}$ recursively (in time), by providing a method of determining $\undt{\mu_{e}}_{0,t}$ in terms of $\undt{\mu_e}_{0,t-1}$. 
Let $K^{{\mu_e},l}_{(t-1,s-1)}$ be the $\olivier{(t-1)} \times \olivier{(s-1)}$ submatrix of $K^{\mu_e,l}$ composed of the rows from times $1$ to $t-1$
and the columns from times $1$ to $s-1$. Let the measure $\undt{\mu_e}^{V_0}_{\,0,t} \in \mathcal{M}(\T_{0,t})$  be given by
\[
\undt{\mu_e}_{\,0,t}^{V_0}(dv)= \mu_I(dv_0) \otimes \mathcal{N}_{t}\left(c^{\mu_e}_{1,t},\sigma^2 \rm{Id}_t + K^{{\mu_e},0}_{(t,t)}\right)dv_{1,t} ,
\]
and $\undt{\mu_e}^{(0,l)}_{\,(0,t),(0,s)} \in \mathcal{M}(\T_{0,t} \times \T_{0,s})$ be given by
\begin{multline*}
\undt{\mu_e}_{\,(0,t),(0,s)}^{(0,l)}(dv^0_{0,t}dv^l_{0,s})= \mu_I(dv_0^0) \otimes \mu_I(dv_0^l) \otimes\\
\mathcal{N}_{t+s}((c^{\mu_e}_{1,t},c^{\mu_e}_{1,s}),\sigma^2 {\rm Id}_{t+s}+K^{{\mu_e},(0,l)}_{(t,s)})dv^0_{1,t} dv^l_{1,s},
\end{multline*}
where
\[
K^{{\mu_e},(0,l)}_{(t,s)}=
\left[
\begin{array}{cc}
 K^{{\mu_e},0}_{(t,t)} & K^{{\mu_e},l}_{(t,s)}\\
{}^\dagger  K^{{\mu_e},l}_{(t,s)} & K^{{\mu_e},0}_{(s,s)}
\end{array}
\right],
\]
\of{and the ${}^\dagger$ sign represents the transpose of a matrix or vector.}

The inductive method for calculating $\undt{\mu_e}$ is outlined in the theorem below.
\begin{theorem}\label{theo:effective}
We may characterise $\undt{\mu_e}$ inductively as follows. Initially $\undt{\mu_e}_{\,0} = \mu_I^{\Z}$. Given that we have a complete characterisation of\\
$\left\lbrace \undt{\mu_e}^{(0,l)}_{\,(0,t-1),(0,t-1)},\undt{\mu_e}_{\,0,t-1}^{V_0} : l\in\Z\right\rbrace, $we may characterise 
$\left\lbrace \undt{\mu_e}_{\,(0,t),(0,t)}^{(0,l)},\undt{\mu_e}_{\,0,t}^{V_0} : l\in\Z\right\rbrace$
according to the following identities. For $s\in [1,t]$,
\begin{equation}\label{eq:mean}
c^{\mu_e}_s = \bar{J} \int_{\R^{s}}f\left(\Psi^{-1}(v)_{s-1}\right)\,\undt{\mu_e}_{\,0,s-1}^{V_0}(dv).
\end{equation}
For $1\leq r,s\leq t$, $K^{\mu_e,k}_{rs}=\theta^2 \delta_k \mone_{T} {}^\dagger  \mone_{T}+
\sum_{l=-\infty}^{\infty} \Lambda(k,l) M^{\mu,l}_{rs}.$ Here, for $p=\text{\rm max}(r-1,s-1)$, 
\begin{equation}\label{eq:Mmue0}
M^{{\mu_e},0}_{rs} = 
\int_{\R^{p+1}} f(\Psi^{-1}(v)_{r-1}) \times f ( \Psi^{-1}(v)_{s-1})\, \undt{\mu_e}_{\,0,p}^{V_0}(dv),
\end{equation}
and for $l \neq 0$
\begin{equation}\label{eq:Mmue}
M^{{\mu_e},l}_{rs} = 
\int_{\R^{r}\times \R^{s}}f(\Psi^{-1}(v^0)_{r-1}) \times f (\Psi^{-1}(v^l)_{s-1})\, \undt{\mu_e}_{\,(0, r-1),(0,s-1)}^{(0,l)}(dv^0dv^l). 
\end{equation}
\end{theorem}
\begin{remark}
From a practical point of view the $t$-dimensional integral in \eqref{eq:mean} and the ${\rm max}(r,s)$-dimensional integral in \eqref{eq:Mmue0} can be reduced by a change of variable to at most two dimensions. Similarly the $r+s$-dimensional integral in \eqref{eq:Mmue}  can be reduced to at most four dimensions.
\end{remark}
\begin{remark}
If we make the biologically realistic assumption that the synaptic weights are not correlated beyond a certain correlation distance $d\geq 0$, $\Lambda(k,l)=0$ if $k$ or $l$ does not belong to $V_d$ it is seen that the matrixes $K^{\mu_e,\,k}$ are 0 as soon as $k \notin V_d$: thus in this case the asymptotic description of the network of neurons is sparse.
\end{remark}
\vspace{-1.0cm}
\section{Convergence results}\label{sec:convergence}\vspace{-0.5cm}
We use the Large Deviation Principle proved in \cite{faugeras-maclaurin:14,faugeras-maclaurin:14d} to establish convergence results for the measures $\Pi^n$, $Q^{V_n}$ and $Q^{V_n}(J^n)$.
\begin{theorem}\label{theo:weakconv}
$\Pi^n$ converges weakly to $\delta_{\mu_e}$, i.e., for all $\Phi \in \mathcal{C}_b(\mathcal{M}_\mS(\T^\Z))$,
\[
\lim_{n \to \infty} \int_{\T^{V_n}} \Phi(\hat{\mu}_n(u))\,Q^{V_n}(du)=\Phi(\mu_e).
\]
Similarly,
\[
\lim_{n \to \infty} \int_{\T^{V_n}} \Phi(\hat{\mu}_n(u))\,Q^{V_n}(J^n)(d\olivier{u})=\Phi(\mu_e)\quad \J \quad \text{almost surely}
\]
\end{theorem}
\begin{pf}
The proof of the first result follows directly from the existence of an LDP for the measure $\Pi^n$, see theorem 3.1 in \cite{faugeras-maclaurin:14}, and is a straightforward adaptation of the one in \cite[Theorem 2.5.1]{moynot:99}. The proof of the second result uses the same method, making use of theorem \ref{theo:asLDP} below.
\qed
\end{pf}\vspace{-0.5cm}
We can in fact obtain the following quenched convergence analogue of the usual lower bound inequality in the definition of a Large Deviation Principle. 
\begin{theorem}\label{theo:asLDP}
For each closed set $F$ of $\mmeas(\T^\Z)$ and for almost all $J$
\[
\lsup{n} \frac{1}{|V_n|} \log \left[Q^{V_n}(J^n) (\hat{\mu}_n \in F)\right]  \leq -\inf_{\mu \in F} H(\mu).
\]
\end{theorem}
\begin{pf}
The proof is a combination of Tchebyshev's inequality and the Borel-Cantelli lemma and is an adaptation of the one in \cite[Theorem 2.5.4, Corollary 2.5.6]{moynot:99}.
\qed
\end{pf}
We define $\check{Q}^{V_n}(J^n) = \frac{1}{|V_n|}\sum_{j \in V_n} Q^{V_n}(J^n)\circ \mS^{-\olivier{j}}$, where we recall the shift operator $\mS$. Clearly $\check{Q}^{V_n}(J^n)$ is in
$\mmeas(\T^{V_n})$. \of{We define $\check{Q}^{V_n}$ to be the expectation of $\check{Q}^{V_n}(J^n)$, with respect to the synaptic weights $J$.}
\begin{corollary}\label{cor:weakconv}
Fix $m$ and let $n> m$. For almost every $J$ and all $h \in \mathcal{C}_b(\T^{V_m})$,
\begin{align*}
\lim_{n \to \infty} \int_{\T^{V_m}}h(u)\,\check{Q}^{V_n,V_m}(J^{n})(du)&=\int_{\T^{V_m}} h(u)\,\mu_e^{V_m}(du).\\
\lim_{n \to \infty} \int_{\T^{V_m}}h(u)\,Q^{V_n,V_m}(du)&=\int_{\T^{V_m}} h(u)\,\mu_e^{V_m}(du).
\end{align*}
That is, the $V_m^{th}$ marginals $\check{Q}^{V_n,V_m}(J^n)$ and $Q^{V_n,V_m}$ \of{ of respectively $\check{Q}^{V_n}(J^N)$ and $\check{Q}^{V_n}$} converge weakly to $\mu_e^{V_m}$ as $n\to\infty$. 
\end{corollary}
\begin{pf}
It is sufficient to apply theorem \ref{theo:weakconv} in the case where $\Phi$ in $\mathcal{C}_b(\mathcal{M}_\mS(\T^\Z))$ is defined by
\[
\Phi(\mu)=\int_{\T^{V_m}} h\,d\mu^{V_m}
\]
and to use the fact that $Q^{V_n},\check{Q}^{V_n}(J) \in \M_\mS(\T^{V_n})$.
\qed
\end{pf}
We now prove the following ergodic-type theorem. We may represent the ambient probability space by $\mathfrak{W}$, where $\omega\in\mathfrak{W}$ is such that $\omega = (J_{ij},B_{t}^j,u^j_0)$, where $i,j\in\Z$ and  $0 \leq t \leq T-1$, recall \eqref{eq:U}. We denote the probability measure governing $\omega$ by $\mathfrak{P}$. Let $u^{(n)}(\omega)\in \T^{V_n}$ be defined by \eqref{eq:U}. As an aside, we may then understand $Q^{V_n}(J^n)$ to be the conditional law of $\mathfrak{P}$ on $u^{(n)}(\omega)$, for given $J^n$. 
\begin{theorem}\label{theo:ergodic}
Fix $m > 0$ and let $h\in C_b(\T^{V_m})$. For $u^{(n)}(\omega)\in\T^{V_n}$ (where $n>m$) $\mathfrak{P}$ almost surely,
\begin{equation}\label{eq:Mthmarginalconverge}
\lim_{n\to\infty} \frac{1}{|V_n|} \sum _{j\in V_n} h\left(\pi^{V_m}(\mS^j u^{(n)}(\omega))\right) = \int_{\T^{V_m}}h(u)d\mu_e^{V_m}(u),
\end{equation}
where $\pi^{V_m}$ is the projection onto $V_m$. Hence $\hat{\mu}_n(u^{\olivier{(n)}}(\omega))$ converges $\mathfrak{P}$-almost-surely to $\mu_e$.
\end{theorem}
\begin{pf}
Our proof is an adaptation of \cite{moynot:99}. We may suppose without loss of generality that $\int_{\T^{V_m}}h(u) d\mu_e^{V_m}(u) = 0$. For $p>1$ let 
\begin{equation*}
F_p = \left\lbrace \mu\in\M_\mS(\T^{\Z})| \left|\int_{\T^{V_m}}h(u)\,\mu^{V_m}(du)\right| \geq \frac{1}{p} \right\rbrace.
\end{equation*}
Since $\mu_e\notin F_p$, but it is the unique zero of $H$, it follows that $\inf_{F_p}H = m > 0$. Thus by theorem 3.1 in \cite{faugeras-maclaurin:14} there exists an $n_0$, such that for all $n>n_0$,
\begin{equation*}
Q^{V_n}\left(\hat{\mu}_n\in F_p\right) \leq \exp\left(-m|V_n|\right). 
\end{equation*}
However
\begin{equation*}
\mathfrak{P}\left(\omega |\hat{\mu}_n(u^{(n)}(\omega))\in F_p\right) = Q^{V_n}\left( u | \hat{\mu}_n(u)\in F_p\right).
\end{equation*}
Thus 
\begin{equation*}
\sum_{n=0}^{\infty}\mathfrak{P}\left(\omega |\hat{\mu}_n(u^{(n)}(\omega))\in F_p\right) < \infty.
\end{equation*}
We may thus conclude from the Borel-Cantelli Lemma that $\mathfrak{P}$ almost surely, for every $\omega\in\mathfrak{W}$, there exists $n_p$ such that for all $n\geq n_p$,
\[
\left|\frac{1}{|V_n|}\sum_{j\in V_n} h\left(\pi^{V_m} (\mS^j u^{(n)}(\omega))\right)\right| \leq \frac{1}{p}. 
\]
This yields \eqref{eq:Mthmarginalconverge} because $p$ is arbitrary. The convergence of $\hat{\mu}_n(u^{\olivier{(n)}}(\omega))$ is a direct consequence of \eqref{eq:Mthmarginalconverge}, since this means that each of the $V_m^{th}$ marginals converge.
\qed
\end{pf}\vspace{-0.5cm}





%
%
%
%

\end{document}